\documentclass[3p]{elsarticle}
\usepackage[table]{xcolor}
\usepackage{natbib}

\usepackage{lineno,hyperref}

\usepackage{latexsym,amsmath,amssymb,amsthm,amsfonts,graphicx,graphics,epsfig,url,color}
\usepackage{tikz}
\usetikzlibrary{positioning}
\usepackage{float}
\usepackage{algorithm}
\makeatletter
\renewcommand{\ALG@name}{Training Algorithm}
\makeatother
\usepackage{algorithmic}
\usepackage{pdflscape}
\usepackage{afterpage}
\modulolinenumbers[5]
\newcommand{\mbf}{\mathbf}

\journal{Physica D}
\makeatletter
\def\ps@pprintTitle{%
 \let\@oddhead\@empty
 \let\@evenhead\@empty
 \def\@oddfoot{\centerline{\thepage}}%
 \let\@evenfoot\@oddfoot}

\begin{document}

\begin{frontmatter}

\title{Supervised Learning Algorithms for Controlling \\Underactuated Dynamical Systems}

\author[mymainaddress]{Bharat Monga\corref{mycorrespondingauthor}}
\ead{monga@ucsb.edu}

\author[mymainaddress,mysecondaryaddress]{Jeff Moehlis}
\cortext[mycorrespondingauthor]{Corresponding author}
\ead{moehlis@engineering.ucsb.edu}

\address[mymainaddress]{Department of Mechanical Engineering, Engineering II Building, University of California Santa Barbara, Santa Barbara, CA 93106, United States}
\address[mysecondaryaddress]{Program in Dynamical Neuroscience, University of California Santa Barbara, Santa Barbara, CA 93106, United States}

\begin{abstract}
Control of underactuated dynamical systems has been studied for decades in robotics, and is now emerging in other fields such as neuroscience. Most of the advances have been in model based control theory, which has limitations when the system under study is very complex and it is not possible to construct a model. This calls for data driven control methods like machine learning, which has spread to many fields in the recent years including control theory. However, the success of such algorithms has been dependent on availability of large datasets. Moreover, due to their black box nature, it is challenging to analyze how such algorithms work, which may be crucial in applications where failure is very costly. In this paper, we develop two related novel supervised learning algorithms. The algorithms are powerful enough to control a wide variety of complex underactuated dynamical systems, and yet have a simple and intelligent structure that allows them to work with a sparse data set even in the presence of noise. Our algorithms output a bang-bang (binary) control input by taking in feedback of the state of the dynamical system. The algorithms learn this control input by maximizing a reward function in both short and long time horizons. We demonstrate the versatility of our algorithms by applying them to a diverse range of applications including: switching between bistable states, changing the phase of an oscillator, desynchronizing a population of synchronized coupled oscillators, and stabilizing an unstable fixed point\footnote{Implementation of the algorithms in this article is available at \url{https://github.com/bharatmonga/Supervised-learning-algorithms}.}. For most of these applications we are able to reason why our algorithms work by using traditional dynamical systems and control theory. We also compare our learning algorithms with some traditional control algorithms, and reason why our algorithms work better. 
\end{abstract}


\end{frontmatter}

Underactuated dynamical systems are systems with fewer actuators/controls than the dimensionality of the state space of the system. Such systems are ubiquitous in a variety of fields including physics, chemistry, biology and engineering. There have been numerous advances made on controlling such systems, with much of the work in robotics \cite{under_mec1,under_mec2,under_mec3}. Control in other applications, especially biology \cite{tutorial,tass2003,circ2,circ1}, is on the rise as it provides promising treatment strategies for several disorders such as Parkinson's disease, cardiac arrhythmias, and jet lag. Most of these control methods, both in robotics and biology, are based on traditional model based control theory and optimal control. 

Such methods work well when it is possible to model the dynamics of the system accurately, which is very difficult as the systems become complicated, especially in neuroscience applications where the dynamics of a single neuron may change rapidly depending on the response from other neurons in the network. Even if an accurate model could be built to describe the dynamics of such a system, developing a classical model based control for such an underactuated system is a challenging task. If the parameters of the system change with time, or if the model doesn't describe the dynamics accurately, the theoretical control guarantees like stability and boundedness may not apply in real systems \cite{HOU20133,anderson2005}. This calls for the development of data driven control algorithms that can learn to control the system without explicitly using a model.

Artificial intelligence algorithms are able to learn to control dynamical systems by using deep neural networks. Such algorithms have been used for a long time \cite{old_neural1,old_neural2}, but with the availability of large data sets, improvement in deep learning architectures, optimization methods, and cheap computation costs, their use is on the rise \cite{Jordan255}. However, the success of such algorithms has been largely dependent on availability of large datasets \cite{nathan}, which can be limited in fields like neuroscience where the cost of obtaining human/animal brain data is very high. Moreover, the black box nature of such methods makes their analysis difficult. Such analysis would be important in tasks where failure is very costly. Another limitation of such methods is their inability to take advantage of the inherent dynamics of the system to achieve the task, which limits their performance.

All these limitations call for a new machine learning control algorithm that doesn't rely on large amounts of data, is easy to understand, and can take advantage of the underlying dynamics in achieving the task. In this article, we have developed two related novel supervised learning algorithms based on these three goals. Our algorithms are powerful enough to control a wide variety of complex underactuated dynamical systems, and yet have a simple structure so one can understand how they work using dynamical systems and control theory foundations. Their simple yet intelligent structure also allows them to effectively achieve the control objective by training on a sparse data set, even in the presence of noise. Our algorithms output a bang-bang (binary) control input by taking in feedback of the state of the dynamical system. The algorithms learn this control input by maximizing a reward function in both short and long time horizons. We demonstrate the versatility of our algorithms by applying them to a diverse range of underactuated dynamical systems including: switching between bistable states, changing the phase of an oscillator, desynchronizing a population of synchronized coupled oscillators, and stabilizing an unstable fixed point of a dynamical system. For most of these applications we are able to reason why our algorithms work by using traditional dynamical systems and control theory. We compare our algorithms with traditional control algorithms and reason why our algorithms work better, especially because they learn to take advantage of the underlying system dynamics in achieving the control objective. We carry out a robustness analysis to demonstrate the effectiveness of our algorithms even in the presence of noise.

This article in organized as follows. In Section \ref{learning_alg}, we develop our supervised learning algorithms to output a binary control. We demonstrate our first supervised learning algorithm by controlling underactuated bistable dynamical systems in Section \ref{bistable}, and compare our algorithm to a fully actuated control. In Section \ref{single_osc}, we illustrate the effectiveness of our second supervised learning algorithm by controlling the phase of a single oscillator and comparing the algorithm to a model based optimal control algorithm. We further demonstrate the versatility of our second supervised learning algorithm by using it to desynchronize a population of synchronized coupled oscillators in Section \ref{coupled}. In Section \ref{unstable}, we apply our second algorithm to stabilize an unstable fixed point of an underactuated dynamical system, and compare the algorithm to a model based control algorithm. To demonstrate the applicability of our algorithms in a more realistic setting, we show how their intelligent structure allows them to perform well in the presence of noise. Section \ref{conclusion} summarizes our work and suggests future extensions and tools. \ref{model_param} lists the mathematical models used in this article. In \ref{phase_red}, we give background on phase reduction relevant for the model based control comparison and validation of our second learning algorithm in Sections \ref{single_osc} and \ref{coupled}, respectively.

\section{Supervised Learning Algorithms}\label{learning_alg}
In this section, we develop our supervised learning algorithms to control a diverse range of underactuated dynamical systems. The development includes two important steps: the first step is generating an appropriate training data set, and the second step is feeding that data set into a binary classifier to control the dynamical system. We call the first step the training algorithm, as it generates the training data by maximizing a reward function. We show in the rest of the paper how choosing an appropriate reward function can be used to take advantage of the inherent dynamics to control a variety of dynamical systems. For the second step, we design a locally weighted binary classifier that takes in the state of the dynamical system as input, and outputs a binary control input in real time. The local nature of our classifier allows it to effectively interpolate between nonlinear boundaries inherent in our data set, and thus plays an important role in the control. The combination of these two steps results in our supervised learning algorithms described below and shown in Figure \ref{flowchart}.

We consider an underactuated dynamical system with an additive control input $\pi\left(\mbf{x}(t)\right)$ as
\begin{eqnarray}
\frac{d }{d t}\mbf{x}(t) = F(\mbf{x}(t)) +[\pi\left(\mbf{x}(t)\right),0_{n-1}]^T, \qquad \mbf{x}(t) \in R^n,\label{udxdt}
\end{eqnarray}
where $0_{n-1}$ is an $n-1$ dimensional zero vector. Thus, the control input depends on the full state of the dynamical system, and only directly affects the first state of the dynamical system. The control input is binary in nature having two values $\{u_1>0,u_2\}$, which can be chosen differently for different applications. For our first algorithm, we take $u_2=0$, and thus the control can be thought of having an ``ON'' state with the value $u_1$ and an ``OFF'' state with value 0.  For our second algorithm, we take the control to be a bang bang control with $u_2=-u_1$. Both algorithms learn the control input $\pi\left(\mbf{x}(t)\right)$ as a function of the state to achieve a particular control objective. They do so by learning from the data generated by sampling a model describing the underlying dynamics of the system. To demonstrate our algorithm in this article, we use an analytical model ($F(\mbf{x}(t))$) which generates our training data. In case a model is not available in an application, one can still use the same algorithms by obtaining training data by direct measurement of the states of the system at different times. Below we describe our supervised learning algorithms in more detail. 

\subsection{Training Algorithm 1}
Our first supervised learning algorithm outputs a ``ON'' and ``OFF'' binary control input. The algorithm learns what control input to output to achieve a certain control objective by maximizing a reward function $\mathcal{R}(\mbf{x}(t))$ which needs to be carefully designed to achieve a control objective in a particular application.

We sample a state $\mbf{x}(0)$ randomly from the state space of the dynamical model of the system, and evolve the state forward in time for short time $\Delta t$ with control state OFF. If $\mathcal{R}(\mbf{x}(\Delta t)) \ge \mathcal{R}(\mbf{x}(0))$, we set the control policy for state $\mbf{x}(0)$ as OFF. If it is not, we again evolve the initial state forward for the same time but with control ON, and compare $\mathcal{R}(\mbf{x}(\Delta t))$ for both control ON and OFF. Whichever control policy maximizes the reward is set for the sampled state $\mbf{x}(0)$. We repeat the process $N$ times by sampling more initial states randomly. The training algorithm is summarized below:

\begin{algorithm}[H]
 \caption{}
  \label{super1}
    \begin{algorithmic}

  \STATE Initialize $\mathcal{X}$ as zeros(N,length($\mbf{x}$)) and $\mathcal{U}$ as zeros(N,1)
  \FOR{i=1,N}
    \STATE Randomly sample \ $\mbf{x}(0)$
    \STATE Compute $\mbf{x}(\Delta t)$ and $\mathcal{R}(\mbf{x}(\Delta t))$ with control OFF
    \IF{$\mathcal{R}(\mbf{x}(\Delta t)) \ge \mathcal{R}(\mbf{x}(0))$} 
         \STATE {Set policy for $\mbf{x}(0)$ as OFF}
     \ELSE
     		\STATE Compute $\mbf{x}(\Delta t)$ and $\mathcal{R}(\mbf{x}(\Delta t))$ with control ON
    
      \IF{$\mathcal{R}(\mbf{x}(\Delta t))$ with control ON $ > \mathcal{R}(\mbf{x}(\Delta t))$ with control OFF}
      		\STATE Set policy for $\mbf{x}(0)$ as ON
      \ELSE
        \STATE Set policy for $\mbf{x}(0)$ as OFF
       \ENDIF
       \ENDIF
      \STATE Assign $\mathcal{X}[i,:]$ as $\mbf{x}(0)$ and $\mathcal{U}[i,:]$ as the policy
    \ENDFOR
    \RETURN $\mathcal{X}$, $\mathcal{U}$
  \end{algorithmic}
\end{algorithm}

Such an algorithm takes advantage of the underlying dynamics by letting the trajectories evolve without any control and only switching ``ON'' the control when necessary. This makes our algorithm highly energy efficient. Such an algorithm is very suitable for controlling bistable dynamical systems where the objective is for the trajectory to converge to a particular stable state of the system, or to switch from one stable state to another. The control can switch OFF when the trajectory enters the region of attraction of the desired stable state, and let the dynamics take the trajectory to the desired state.

\subsection{Training Algorithm 2}
Our second supervised learning algorithm outputs a bang-bang control input which can be used to control a variety of dynamical systems, including coupled oscillators. The algorithm learns what control input to output to achieve a certain control objective by maximizing a reward function $\mathcal{R}(\mbf{x}(\Delta t))$, which needs to be carefully designed to achieve a control objective in a particular application.

We sample a state $\mbf{x}(0)$ randomly from the state space of the dynamical model of the system, and evolve the state forward in time for short time $\Delta t$ with both control $u_1$ and $-u_1$. In both scenarios we let the state evolve further in time with zero control until some event occurs and measure the timing of this event. The reward $\mathcal{R}(\mbf{x}(\Delta t))$ is dependent on the timing of this event.  Whichever policy ($u_1$ or $-u_1$) maximizes this reward is set for that sampled state. We repeat the process $N$ times by sampling more states randomly. The algorithm is summarized below:

\begin{algorithm}[H]
 \caption{}
  \label{super2}
  \begin{algorithmic}
  \STATE Initialize $\mathcal{X}$ as zeros(N,length($\mbf{x}$)) and $\mathcal{U}$ as zeros(N,1)
  \FOR{i=1,N}
    \STATE Randomly sample \ $\mbf{x}(0)$
    \STATE Compute $\mbf{x}(\Delta t)$ and $\mathcal{R}(\mbf{x}(\Delta t))$ with control $u_1$
     \STATE Compute $\mbf{x}(\Delta t)$ and $\mathcal{R}(\mbf{x}(\Delta t))$ with control -$u_1$
      \IF{$\mathcal{R}(\mbf{x}(\Delta t))$ with control $u_1 \ge \mathcal{R}(\mbf{x}(\Delta t))$ with control -$u_1$}
      		\STATE Set policy for $\mbf{x}(0)$ as $u_1$
      \ELSE
        \STATE Set policy for $\mbf{x}(0)$ as -$u_1$   
     \ENDIF
     \STATE Assign $\mathcal{X}[i,:]$ as $\mbf{x}(0)$ and $\mathcal{U}[i,:]$ as the policy
    \ENDFOR
    \RETURN $\mathcal{X}$, $\mathcal{U}$
  \end{algorithmic}
\end{algorithm}

Such control is useful when the objective is to converge to an unstable state of the system, because the control needs to stay ``ON'' (be non-zero) the whole time even when the control objective has been realized, since the trajectory will go back to the stable state otherwise. Here as well the underlying dynamics of the system play a role in our learning algorithm: to determine the control input for a particular initial state, we let the dynamics evolve the trajectory until an event occurs.

Both of these training algorithms generate data comprising a set of $N$ sampled states of the dynamical system $\mathcal{X}$, and a set of the corresponding control inputs $\mathcal{U}$. However, we need to know the control input $\pi\left(\mbf{x}(t)\right)$ as a function of a general trajectory $\mbf{x}(t)$ of the system, since the trajectory is not restricted to these sampled states. This is achieved with our binary classifier, that takes an input as the state of the dynamical system $\mbf{x}(t)$ and outputs the corresponding control policy $\pi\left(\mbf{x}(t)\right)$ based on this generated data.

\subsection{Binary Classifier}
Our training algorithms generate data comprising a set $\mathcal{X}$ of $N$ sampled states of the dynamical system, and a set of the corresponding control inputs $\mathcal{U}$. Based on this information we build a locally weighted binary classifier that takes as input as the instantaneous state of the dynamical system $\mbf{x}(t)$, and outputs the corresponding control policy $\pi\left(\mbf{x}(t)\right)$ to be applied at that instant. 

We assign each element of the set $\mathcal{X}$ with a weight 
\begin{eqnarray*}
w_i (\mbf{x}(t))= \exp\left(-\frac{\left|\mbf{x}(t)-\mathcal{X}_i\right|^2}{2\tau}\right),\quad i=1,2,\ldots, N,
\end{eqnarray*}
where $\mathcal{X}_i$ represents the $i^{th}$ sampled state stored in the set $\mathcal{X}$. Thus a sampled state is given a higher weight if its closer to $\mbf{x}(t)$, and a lower weight if it is further away from $\mbf{x}(t)$. These weights are normalized so that 
\begin{eqnarray*}
\sum_{i=1}^{N} w_i (\mbf{x}(t))= 1.
\end{eqnarray*}
For the first algorithm, the classifier outputs $\pi\left(\mbf{x}(t)\right)=u_1$ (``ON'') if
\begin{eqnarray}
\sum_{i=1}^{N} w_i (\mbf{x}(t))\mathcal{U}_i>0.5u_1,\label{bc1}
\end{eqnarray}
and $\pi\left(\mbf{x}(t)\right)=0$ (``OFF'')  otherwise. Similarly for the second algorithm, the classifier outputs $\pi\left(\mbf{x}(t)\right)=u_1$  if
\begin{eqnarray}
\sum_{i=1}^{N} w_i (\mbf{x}(t)) \mathcal{U}_i>0,\label{bc2}
\end{eqnarray}
and $\pi\left(\mbf{x}(t)\right)=-u_1$ otherwise. Here $\tau$ is a hyperparameter, which can be thought of as a bandwidth parameter. Smaller $\tau$ decreases the influence of faraway data points in local interpolation. Choosing a very small value can lead to overfitting and tight boundaries around the sampled points. On the other hand, a large value can bias the decision boundary in favor of a larger cluster of data points. We run experiments with different $\tau$ values and choose the one which enables our classifier to label the training data correctly, and also produces a smooth decision boundary between the clusters. The critical number on the right hand side of equations (\ref{bc1}) and (\ref{bc2}) can be thought of as another hyperparameter which would influence the decision boundary together with the bandwidth  parameter. We keep this parameter fixed for all our case studies, and find that a value of $0.5u_1$ and 0 respectively prove effective in outputting a smooth boundary decision around the cluster of data points for the examples that we consider. Choosing a larger (resp., smaller) critical value would extend the boundary in favor of clusters with labels OFF / negative (resp., ON / positive) control. We do not claim that the chosen hyperparameter values are optimal, but we believe that they are close to optimal for our examples, as the binary classifier is able to correctly label the training data and produce a smooth boundary, and thus is able to achieve the control objectives demonstrated in later sections of the article.

\begin{figure}[H]
\begin{center}
\includegraphics[width=0.5\textwidth]{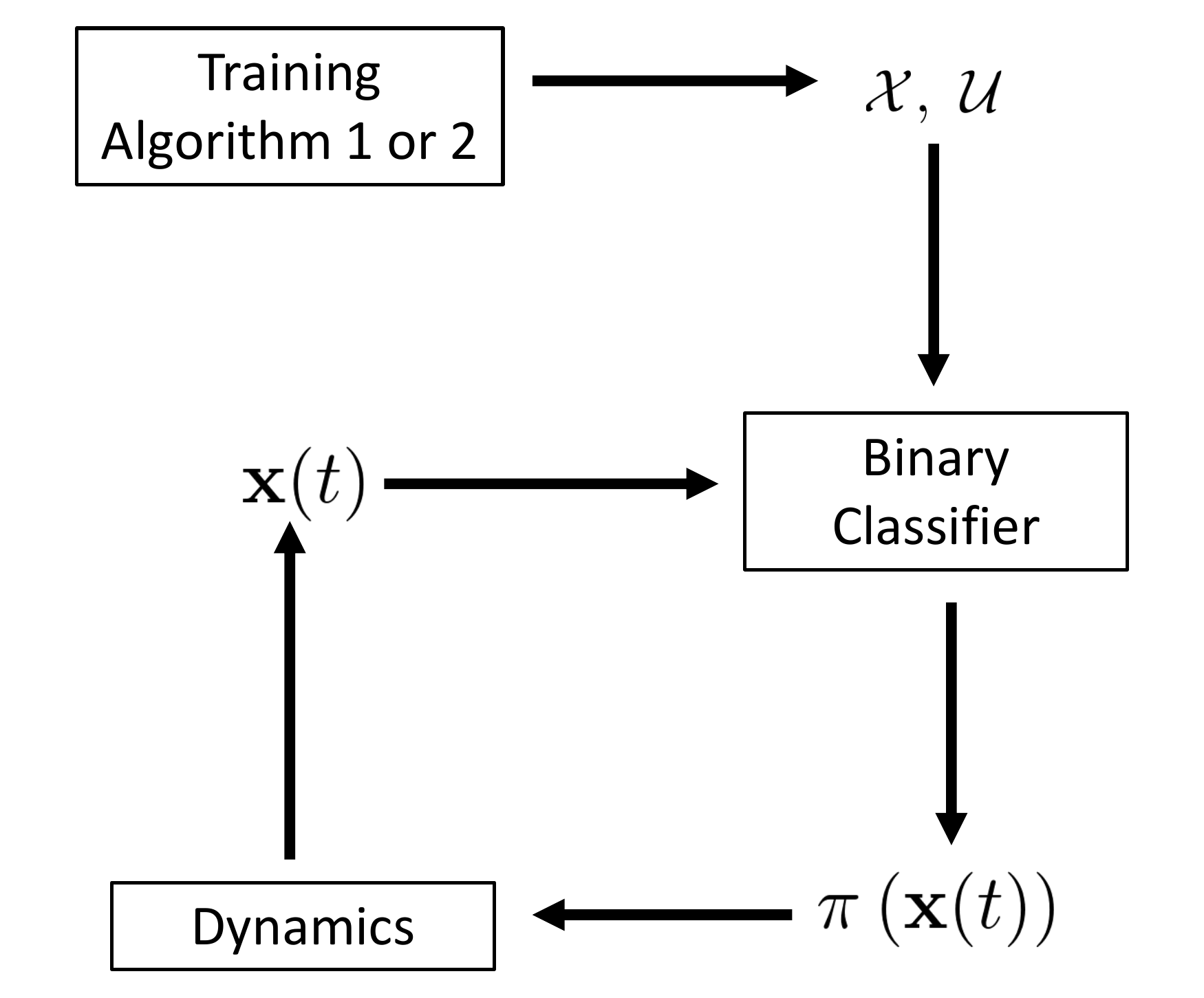}
\end{center}
\caption{Flowchart of the two supervised learning algorithms.}
\label{flowchart}
\end{figure}

With $\pi\left(\mbf{x}(t)\right)$ defined to be the output of this binary classifier, we simulate the dynamical system from (\ref{udxdt}) starting from a random initial condition and find that our supervised learning algorithms are able to achieve the desired control objectives, while simultaneously maximizing the designed reward functions. The entire algorithm is depicted in the flowchart in Figure \ref{flowchart}.

Note that in our algorithm, the sets $\mathcal{X}, \ \mathcal{U}$ need to be computed only once for a given dynamical system, whereas the control input $\pi\left(\mbf{x}(t)\right)$ is computed by the binary classifier at every timestep.

\section{Bistable Dynamical Systems}\label{bistable}
Bistability is widely found in neural systems \cite{neural} and cardiac arrhythmia \cite{cardiac}, and is used in digital electronics for storing binary information, in mechanical switches for transitioning between ON and OFF states, and in multivibrators, Schmitt trigger circuits, and even optical systems \cite{optical}. It is the key mechanism for understanding several cellular processes including those associated with the onset and treatment of cancer \cite{cancer}. In this section we apply our first supervised learning algorithm to control underactuated bistable dynamical systems. The control objective is for the trajectory to converge to a particular stable fixed point of the system starting anywhere in the state space (including in the basin of attraction of the other stable state). Such a control objective is relevant for several applications such as biocomputing, gene therapy, and treatment of cancer \cite{gene,cancer1}, among others.

\subsection{Duffing Oscillator}\label{duff_no_noise}
With the Duffing oscillator \cite{wiggins,duff}, we consider the class of bistable dynamical systems having two stable fixed points ($\mbf{x}^{s1}, \ \mbf{x}^{s2}$), and an unstable fixed point ($\mbf{x}^{u}$). The control objective is for the trajectory to converge to $\mbf{x}^{s2}$ starting anywhere in the state space. The Duffing oscillator is given as:
\begin{eqnarray*}
\dot x&=&y +\pi\left(\mbf{x}(t)\right) ,\\
\dot y &=& x - x^3 -\delta y.
\end{eqnarray*}
For $\delta>0$, the system has two stable fixed points $\mbf{x}^{s1}=(-1,0)$ and $\mbf{x}^{s2}=(1,0)$, and an unstable fixed point $\mbf{x}^{u}=(0,0)$, all shown in Figure \ref{duff}. We take $\delta=0.1$ in our simulations.

\subsubsection{Learning Algorithm}
We choose our reward function to be the negative of the Euclidean distance between the current state and the desired state:
\begin{equation}
\mathcal{R}(\mbf{x}(t)) = -||\mbf{x}(t)-\mbf{x}^{s2}||.\label{reward_bistable}
\end{equation}
Thus the control will make the trajectory converge towards the desired fixed point while increasing the reward to 0. To converge to $\mbf{x}^{s2}$ starting anywhere in the state space, we use our learning algorithm to generate a control policy. The ON (resp., OFF) state of the control policy corresponds to a value of $u_1=4$ (resp., 0). We randomly sample $N=50$ points for generating the sets $\mathcal{X}, \ \mathcal{U}$, and choose $\Delta t = 0.001$ and $\tau=0.4$.

\subsubsection{Results}
\begin{figure}[!t]
\begin{center}
\includegraphics[width=\textwidth]{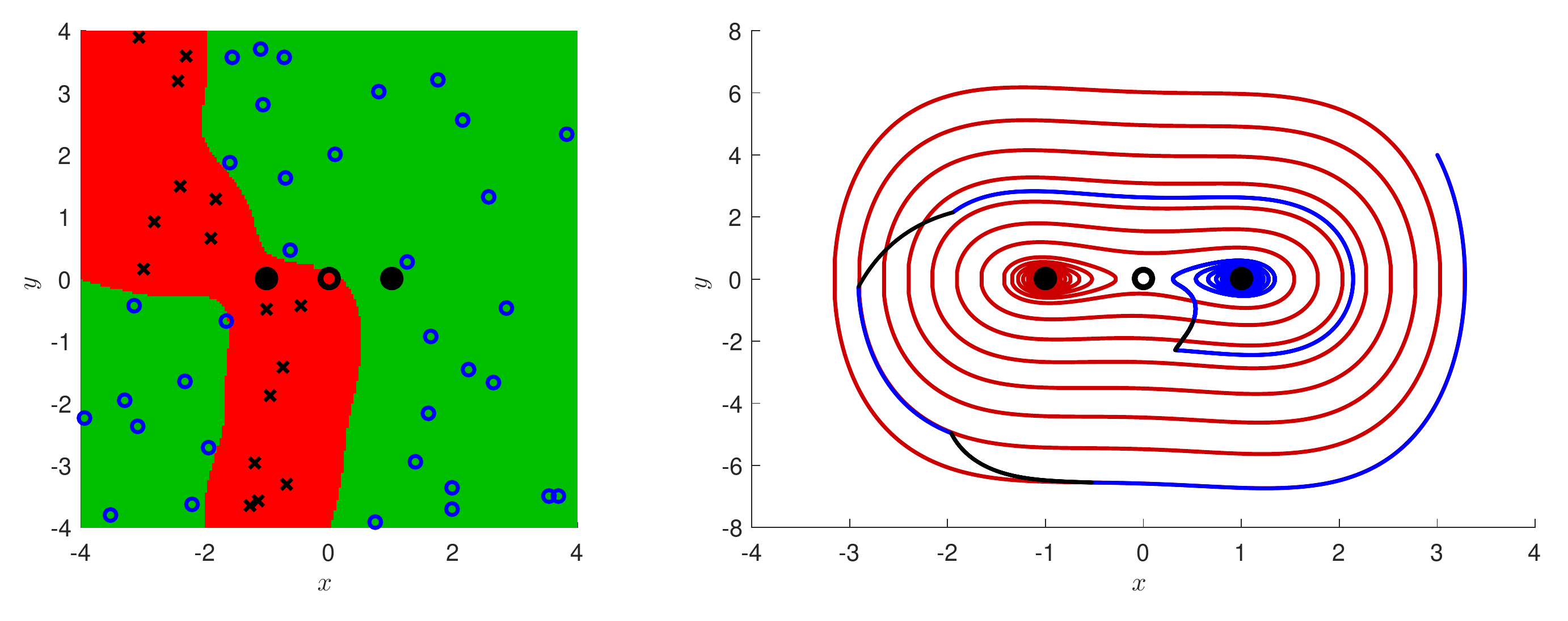}
\end{center}
\caption{Duffing Oscillator ($\delta=0.1$): Solid (resp., open) black circles represent $\mbf{x}^{s1}$, $\mbf{x}^{s2}$, (resp., $\mbf{x}^{u}$). In the left panel, open blue circles (resp., black $\times$'s) represent elements of the set $\mathcal{X}$ where the control policy given by elements of the set $\mathcal{U}$ is OFF (resp., ON). The green (resp., red) region is where the output $\pi\left(\mbf{x}(t)\right)$ of the binary classifier is OFF (resp., ON). In the right panel, the trajectory starts in the region of attraction of $\mbf{x}^{s1}$, and converges to $\mbf{x}^{s2}$ (resp., $\mbf{x}^{s1}$) with (resp., without) control. When the control policy is ON (resp., OFF), the trajectory is plotted in black (resp., blue). The uncontrolled trajectory is plotted in red.}
\label{duff}
\end{figure}

The generated control policy is shown in the left panel of Figure \ref{duff}. Blue open circles (resp., black $\times$'s) represent elements of the set $\mathcal{X}$ where the control policy given by elements of the set $\mathcal{U}$ is OFF (resp., ON). The green (resp., red) region is where output $\pi\left(\mbf{x}(t)\right)$ of the binary classifier is OFF (resp., ON). A controlled and an uncontrolled trajectory starting from same $\mbf{x}(0)$ are shown in the right panel of Figure \ref{duff}. As can be seen in this figure, the control algorithm gradually converges the trajectory to $\mbf{x}^{s2}$ by turning the control ON a few times, whereas the uncontrolled trajectory converges to $\mbf{x}^{s1}$. 

The algorithm generates an energy efficient control policy as the policy is OFF 60.41\% of the total time it takes to drive the trajectory within a ball of radius of 0.45 in the region of attraction of $\mbf{x}^{s2}$. We investigate the robustness of our learning algorithm by testing it on 1000 randomly generated initial conditions, and in all the 1000 cases, the control algorithm is able to converge the trajectories to $\mbf{x}^{s2}$, achieving 100\% accuracy. Choosing $N$ is the crucial task in our learning algorithm. We start with a small $N$ and keep increasing it until the algorithm achieves 100\% effectiveness. $N=50$ points turns out to be appropriate for the Duffing oscillator as choosing a lower number of points leads to underfitting, and choosing a higher number of points leads to a higher computational cost.

\subsection{Reduced Hodgkin-Huxley model}\label{hh}
\begin{figure}[!t]
\begin{center}
\includegraphics[width=\textwidth]{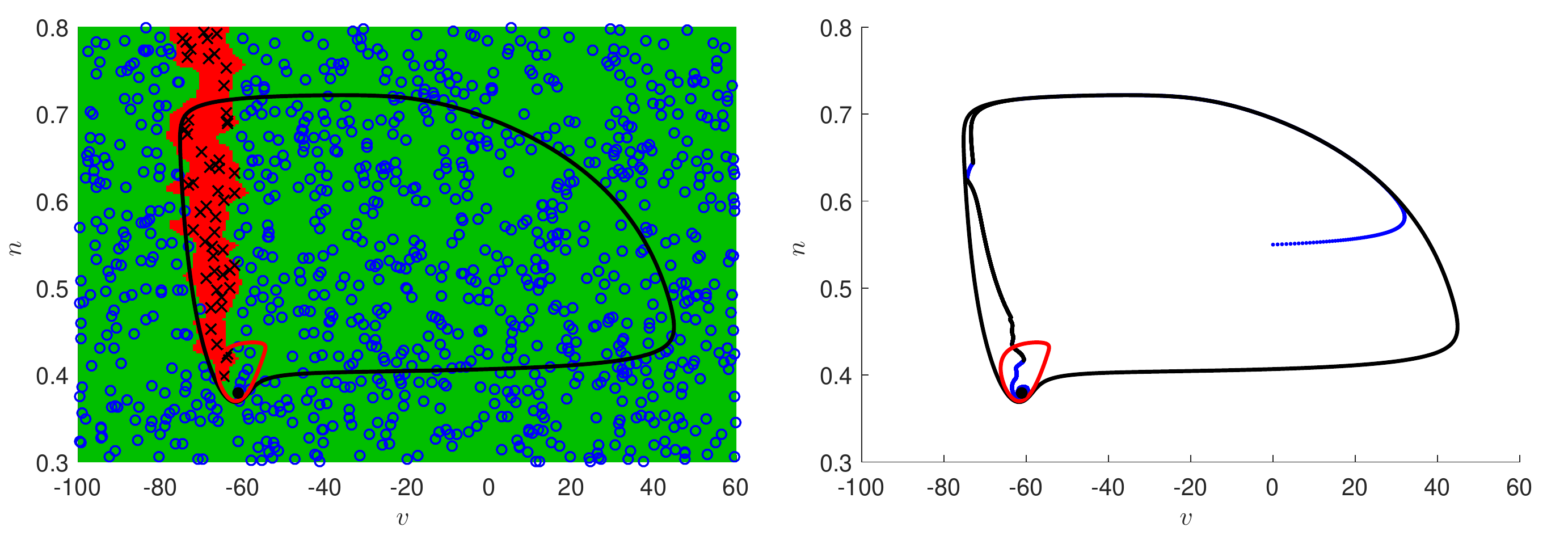}
\end{center}
\caption{Reduced Hodgkin-Huxley model: The black and red curves are $\mbf{x}^{s1}(t)$ and $\mbf{x}^{u}(t)$, respectively. The black point in the bottom left corner of figure panels is $\mbf{x}^{s2}$. In the left panel, small black circles (resp., black $\times$'s) represent elements of the set $\mathcal{X}$ where the control policy given by elements of the set $\mathcal{U}$ is OFF (resp., ON). Green (resp., red) regions are where the output $\pi\left(\mbf{x}(t)\right)$ of the binary classifier is OFF (resp., ON). In the right panel, the trajectory starts outside $\mbf{x}^{u}(t)$, and converges to $\mbf{x}^{s2}$. When the control policy is ON (resp., OFF), the trajectory is plotted in black (resp., blue) color.}
\label{hh_results}
\end{figure}

With the reduced Hodgkin-Huxley model \cite{canard,Keener2009,huxley}, we consider the class of bistable dynamical systems having a stable periodic orbit $\mbf{x}^{s1}(t)$, an unstable periodic orbit $\mbf{x}^{u}(t)$, and a stable fixed point $\mbf{x}^{s2}$. The model is given as
\begin{eqnarray*}
\dot v&=& \left(I - g_{Na}(m_\infty(v))^3(0.8-n)(v-v_{Na}) - g_Kn^4(v-v_K)-g_L(v-v_L)\right)/c+\pi\left(\mbf{x}(t)\right),\\ 
\dot n &=& a_n(v)(1-n)-b_n(v) n,
\end{eqnarray*}
where $v$ is the trans-membrane voltage, and $n$ is the gating variable. $I$ is the baseline current, which we take as 6.69 $\mu A/cm^2$, and $\pi\left(\mbf{x}(t)\right)$ represents the applied control policy. For the rest of the parameters, see \ref{hh_param}. In the absence of control input, the system is bistable having $\mbf{x}^{s1}(t)$ with period 14.91 ms, $\mbf{x}^{u}(t)$ with period 14.33 ms, and $\mbf{x}^{s2}=(-61.04,0.38)$, all shown in Figure \ref{hh_results}. 

\subsubsection{Learning Algorithm}\label{hh_la}
The control objective is for the trajectory to converge to the stable fixed point starting anywhere in the state space. Here as well we choose the reward function (\ref{reward_bistable}). Without any control input, a trajectory starting outside $\mbf{x}^{u}(t)$ will converge to $\mbf{x}^{s1}(t)$, and a trajectory starting inside $\mbf{x}^{u}(t)$ will converge to $\mbf{x}^{s2}$. To converge to the stable fixed point starting anywhere in the state space, we use our learning algorithm to generate a control policy. The ON (resp., OFF) state of the control policy corresponds to a value of $u_c=15$ (resp., 0). We sample $N=1000$ points for generating the sets $\mathcal{X}, \ \mathcal{U}$, and choose $\Delta t = 0.001$, and $\tau=0.001$. Because the two state variables $v, \ n$ scale differently, it is important to normalize them for calculating the reward function. This is also important for the binary classifier to work effectively, since it is based on the Euclidean norm. To do this, we subtract from each element of the set $\mathcal{X}$ the mean of the set and then divide each element by the variance of the set. We subtract the same mean and divide by the same variance from the state $\mbf{x}(t)$ that goes in calculating the reward function and also the binary classifier. 
\subsubsection{Results}
The generated control policy is shown in the left panel of Figure \ref{hh_results}. As shown in this figure, the learning algorithm indicates that is it is better to have control ON in the left part of the state space in order to maximize the reward function. In all other regions, the learning algorithm indicates that the control policy should be OFF. Since the control policy is ON in only a small region of the state space, we need to sample 1000 points to accurately determine this region. A controlled trajectory using this policy is shown in the right panel of Figure \ref{hh_results}. The learning algorithm is able to converge the trajectory to the stable fixed point $\mbf{x}^{s2}$ by bypassing the unstable periodic orbit $\mbf{x}^{u}(t)$. The algorithm generates an energy efficient control policy as the policy is OFF 23.81\% of the time it takes for the algorithm to drive the trajectory inside $\mbf{x}^{u}(t)$. We investigate the robustness of our learning algorithm by testing it on 1000 randomly generated initial conditions, and in all the 1000 cases, the algorithm is able to converge the trajectories to $\mbf{x}^{s2}$, achieving 100\% effectiveness. Note that the learning algorithm has no information about the periodic orbits and fixed points of the system; it only works to maximize the reward function. 

\subsubsection{Comparison with fully actuated control}
To further demonstrate energy efficiency of our learning algorithm, we compare it with a fully actuated feedback control given as
\begin{eqnarray}
\frac{d }{d t}\mbf{x}(t) &=& F(\mbf{x}(t)) +U(\mbf{x}(t)), \qquad \mbf{x}(t) \in R^n,\\
U(\mbf{x}(t)) &=& -F(\mbf{x}(t)) -0.2 \left( \mbf{x}(t) - \mbf{x}^{s2} \right),
\end{eqnarray}
which also converges the trajectory to the stable fixed point in the same time frame as our learning algorithm. However the energy required by this algorithm calculated as $\int_0^{t} ||U(\mbf{x}(t))||_2 ^2dt$ comes out to be more than 3 orders of magnitude larger compared to the energy taken by the control obtained from our learning algorithm. This is because our learning algorithm takes advantage of the natural dynamics of the system to drive the trajectory close to the desired point, and turns the control ON only for a short amount of time when its really needed. In contrast, the feedback based control is ON the whole time, even when the trajectory reaches inside the unstable periodic orbit.

\section{Phase Control of an Oscillator}\label{single_osc}
In this section, we use our second algorithm to control a class of underactuated dynamical systems having a stable limit cycle solution $\mbf{x}^{s}(t)$. We seek to maximally increase or decrease the phase of the limit cycle solution by using a bang-bang type control input. The motivation behind such a control objective comes from controlling neurons, where one might want a neuron to spike as quickly as possible subject to a constraint on the magnitude of the allowed input current; this constraint can be due to hardware limitations and/or concern that large inputs might cause tissue damage. Thus, instead of thinking in terms of maximally increasing the phase, one can instead think in terms of maximally decreasing the neuron's spike time.




\subsection{Model}
\begin{figure}[!t]
\begin{center}
\includegraphics[width=\textwidth]{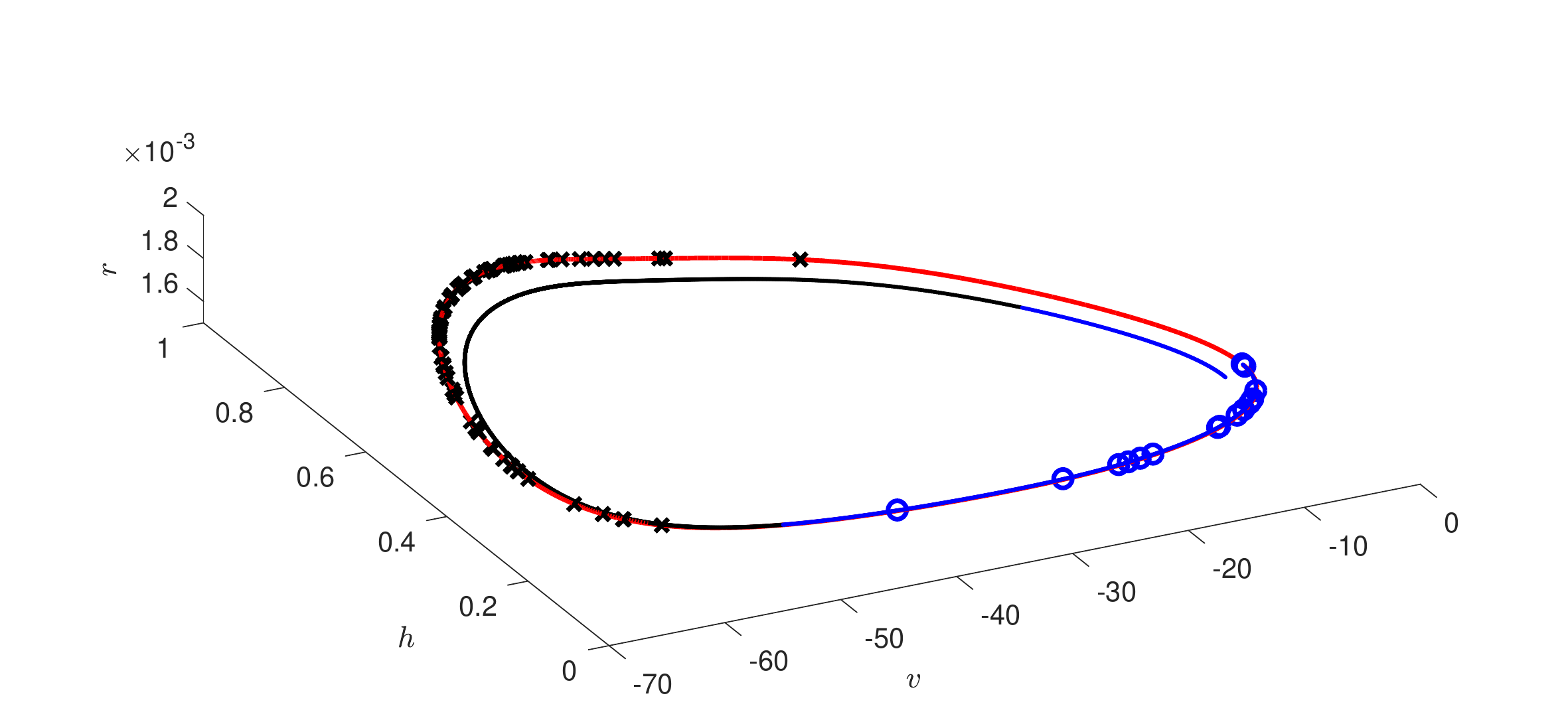}
\end{center}
\caption{Thalamic neuron model: Red curve is the stable limit cycle. Small blue circles (resp., black $\times$'s) represents elements of set $\mathcal{X}$ where control policy given by elements of set $\mathcal{U}$ is $-u_1$ (resp., $u_1$). The controlled trajectory is plotted in blue (where $\pi\left(\mbf{x}(t)\right)=-u_1$) and black (where $\pi\left(\mbf{x}(t)\right)=u_1$).}
\label{th}
\end{figure}

To demonstrate our algorithm, we consider the 3-dimensional thalamic neuron model \cite{rubi04} for the oscillatory behavior of neurons in the thalamus:

\begin{eqnarray}
\dot v&=&\frac{-I_L(v)-I_{Na}(v,h)-I_K(v,h)-I_T(v,r)+I_b }{C_m}+\pi\left(\mbf{x}(t)\right),\label{th1}\\
\dot h&=&\frac{h_{\infty}(v)-h}{\tau_h(v)},\\
\dot r&=&\frac{r_{\infty}(v)-r}{\tau_r (v)},\label{th2}
\end{eqnarray}
where the state $\mbf{x}(t)$ is the tuple $(v,h,r)$, $v$ is the transmembrane voltage, and $h,\ r$ are the gating variables of the neuron. $\pi\left(\mbf{x}(t)\right)$ represents the applied current as the control input. For details of the rest of the parameters, see \ref{thalam_param}. With no control input, these parameters give a stable limit cycle $\mbf{x}^{s}(t)$ with period $T=8.40~ms$ shown in red in Figure \ref{th}.

\subsection{Learning Algorithm}
Here the control objective is to maximally decrease the spike time of the neuron, meaning we want the oscillation to end sooner than it naturally would. We set the reward function as the negative of the neuron's next spike time (the time when the transmembrane voltage $v(t)$ reaches a maximum):
\begin{equation}
\mathcal{R}(\mbf{x}(t)) = -t_{spike}.\label{reward_phase}
\end{equation}
We sample 100 states randomly along the limit cycle and evolve them with both positive and negative control inputs for time $\Delta t=0.001$, and then evolve them further with zero control input until the neuron spikes. Whichever control input attains the minimal $t_{spike}$ (maximizes the reward function) is selected as control policy for that sampled state. We choose $\tau=0.01$. Because the state variables $v, \ h, \ r$ have different dynamic ranges, we normalize the set $\mathcal{X}$ and the state at every time step similar to in Section \ref{hh_la}.

\subsection{Results}
The generated control policy, along with the controlled trajectory, are shown in Figure \ref{th}. As seen in this figure, most of the sampled states need to have a positive control in order to maximize the reward function. This is evident from the left panel of Figure \ref{compare} which plots the corresponding control input. Because of the control, the neuron spikes ($v$ reaches its maximum) in $t_{spike}=7.49 ms$ which is $10.82\%$ decrease in its natural spike time of $8.40~ms$. Thus our algorithm is able to achieve the control objective while keeping the controlled trajectory close to the stable limit cycle solution (see Figure \ref{th}).

\subsection{Model based control comparison}
The dynamics of neural oscillations are highly nonlinear and high dimensional, which makes a model based control formulation very challenging. Phase reduction, a model reduction technique valid close to the limit cycle (see \ref{phase_red} for more details), can overcome these challenges. The neuron spike time control problem was solved as an optimal control problem in \cite{nabi12,tutorial} using phase reduction, which also resulted in a bang bang control with control input given as
\begin{eqnarray}
\pi\left(\mbf{x}(t)\right) &=& -\text{sign} [\mathcal{Z}(\theta)]u_1  \quad \text{for decreasing} \ t_{spike},
\end{eqnarray}
where $u_1$ is the bound chosen by the user and $\mathcal{Z}(\theta)$ is the phase response curve (see \ref{phase_red} for more details) which is a periodic function of $\theta$. Such a control works well, except when the bound $u_1$ is large, where the controlled trajectory can diverge far away from the limit cycle, decreasing the accuracy of phase reduction and making the control based on phase reduction ineffective. Effectiveness of such a control also relies heavily on accurate measurement of the phase response curve, which may not be possible.

\begin{figure}[!t]
\begin{center}
\includegraphics[width=\textwidth]{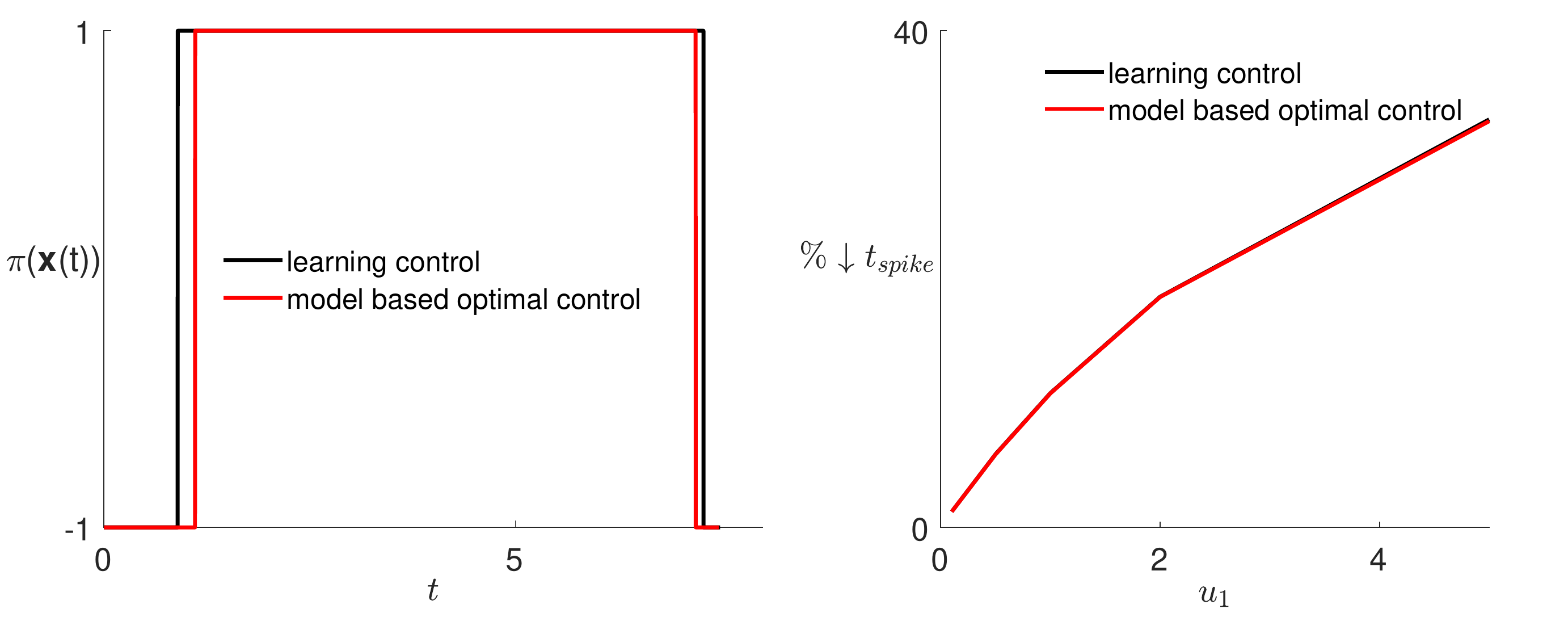}
\end{center}
\caption{Thalamic neuron model: The left panel plots the control input $\pi\left(\mbf{x}(t)\right)$ for our learning algorithm and optimal control algorithm for $u_1=1$. The right panel shows the \% decrease in $t_{spike}$ as a function of $u_1$. }
\label{compare}
\end{figure}

We find that our supervised learning based control outputs a control input very similar to the above model based control, both shown in the left panel of Figure \ref{compare} for $u_1=1$.
We compute $t_{spike}$ as a function of the bound $u_1$ and find that our learning based algorithm does slightly better than the model based algorithm in decreasing $t_{spike}$ (shown in the right panel of Figure \ref{compare}). Both controls are able to decrease $t_{spike}$ more as $u_1$ increases.

\section{Desynchronization of a Population of Coupled Oscillators}\label{coupled}
Populations of coupled oscillators are ubiquitous in applications from physics, chemistry, biology, and engineering \cite{physics,kuramoto,winfree,Keener2009}. The collective behavior of such oscillators varies, and includes synchronization, desynchronization, and clustering in various scenarios. Pathological synchronization of neural oscillations in the thalamus and the subthalamic nucleus (STN) brain region is hypothesized to be one of the causes of motor symptoms for essential and parkinsonian tremor, respectively~\cite{thalamus_syncro,STN_syncro}. Deep brain stimulation (DBS), an FDA approved treatment, has proven to alleviate these symptoms~\cite{DBS_STN,DBS_thalamus} by stimulating the thalamus or the STN brain regions with a high frequency, (relatively) high energy pulsatile waveform, which has been hypothesized to desynchronize the synchronized neurons; see, e.g., \cite{wilson11,clustered}. This has motivated researchers to come up with efficient model dependent control techniques~\cite{tass2003,neural_control2,dan14,monga18,mongaD} to desynchronize the neural oscillations, but also consume less energy, thus prolonging the battery life of the stimulator and preventing tissue damage or side effects caused by the pulsatile stimuli. 

\subsection{Model}
Inspired by such treatment of parkinsonian and essential tremor, we employ our algorithm to desynchronize an initially synchronized population of $M$ coupled thalamic neural oscillations. We consider the 3-dimensional thalamic neuron model \cite{rubi04} for each individual oscillator with added all-to-all electrotonic coupling:

\begin{eqnarray}
\dot v_i&=&\frac{-I_L(v_i)-I_{Na}(v_i,h_i)-I_K(v_i,h_i)-I_T(v_i,r_i)+I_b+\frac{1}{N}\sum_{j=1}^{M}\alpha_{ij}(v_j-v_i) }{C_m}\nonumber\\ &&\qquad  +\pi\left(\mbf{x}(t)\right),\\
\dot h_i&=&\frac{h_{\infty}(v_i)-h_i}{\tau_h(v_i)},\\
\dot r_i&=&\frac{r_{\infty}(v_i)-r_i}{\tau_r (v_i)},
\end{eqnarray}
where $\mbf{x}(t)$ represents the full state (3 $\times M$ dimensional) of the oscillator population. Here, $i = 1,\cdots, M$, where $M$ is the total number of oscillators in the neuron population. $v_i$ is the transmembrane voltage, and $h_i,\ r_i$ are the gating variables of the $i^{th}$ neural oscillator. $\alpha_{ij}$ is the coupling strength between oscillators $i$ and $j$, which are assumed to be electrotonically coupled \cite{coupling} with $\alpha_{ij}=0.01$ for $j\ne i$ and $\alpha_{ii} = 0$ for all $i$. $\pi\left(\mbf{x}(t)\right)$ represents the applied current as the control policy. For details of the rest of the parameters, see \ref{thalam_param}. Note that the same control input $\pi\left(\mbf{x}(t)\right)$ is applied to all of the oscillators. With no control input, these parameters give a synchronized oscillator population with period $T=8.40~ms$. 

\subsection{Learning Algorithm}
We index the individual neural oscillators in the order in which they spike, thus neuron $1$ spikes first and neuron $M$ spikes last. We set the reward function as the absolute value of spike time difference of neuron 1 and $M$: 
\begin{equation}
\mathcal{R}(\mbf{x}(t)) = |t_{spike1}-t_{spikeM}|.\label{reward_osc}
\end{equation}
Since the oscillator population is initially synchronized, this reward is initially a small positive number as all neurons spike very close to each other. We aim to desynchronize the population by maximally increasing this reward function. We consider $M=51$ oscillators in the synchronized population and sample 51 states along the synchronized oscillation. Since the state of the oscillator population is very large ($3\times M$), we take the mean across the population to reduce the dimension of our set $\mathcal{X}$. The $i^{th}$ element of the set $\mathcal{X}$ is given as
\begin{eqnarray}
\mathcal{X}_i=\frac{\sum_{j=1}^{M}\left(v_j,h_j,r_j\right)}{M}.
\end{eqnarray}
We evolve the oscillator population with both positive and negative control inputs for time $\Delta t=0.001$, and evolve them further with zero control input until all neurons in the population spike. Whichever policy attains the maximum reward function is selected for that sampled state $\mathcal{X}_i$. 

The binary classifier takes as input the full high dimensional state of the oscillator population. It then computes the mean of the state across the population and compares it with the sampled mean states to output a control policy $\pi\left(\mbf{x}(t)\right)$. Because the mean of the states $v_j, \ h_j, \ r_j$ scale differently, it is important to normalize them for the binary classifier to work effectively, since it is based on the Euclidean norm. Thus we normalize the set $\mathcal{X}$ and the mean state at every time step similar to in Section \ref{hh_la}. We choose $\tau=0.01$.

\subsection{Results}
\begin{figure}[!t]
\begin{center}
\includegraphics[width=\textwidth]{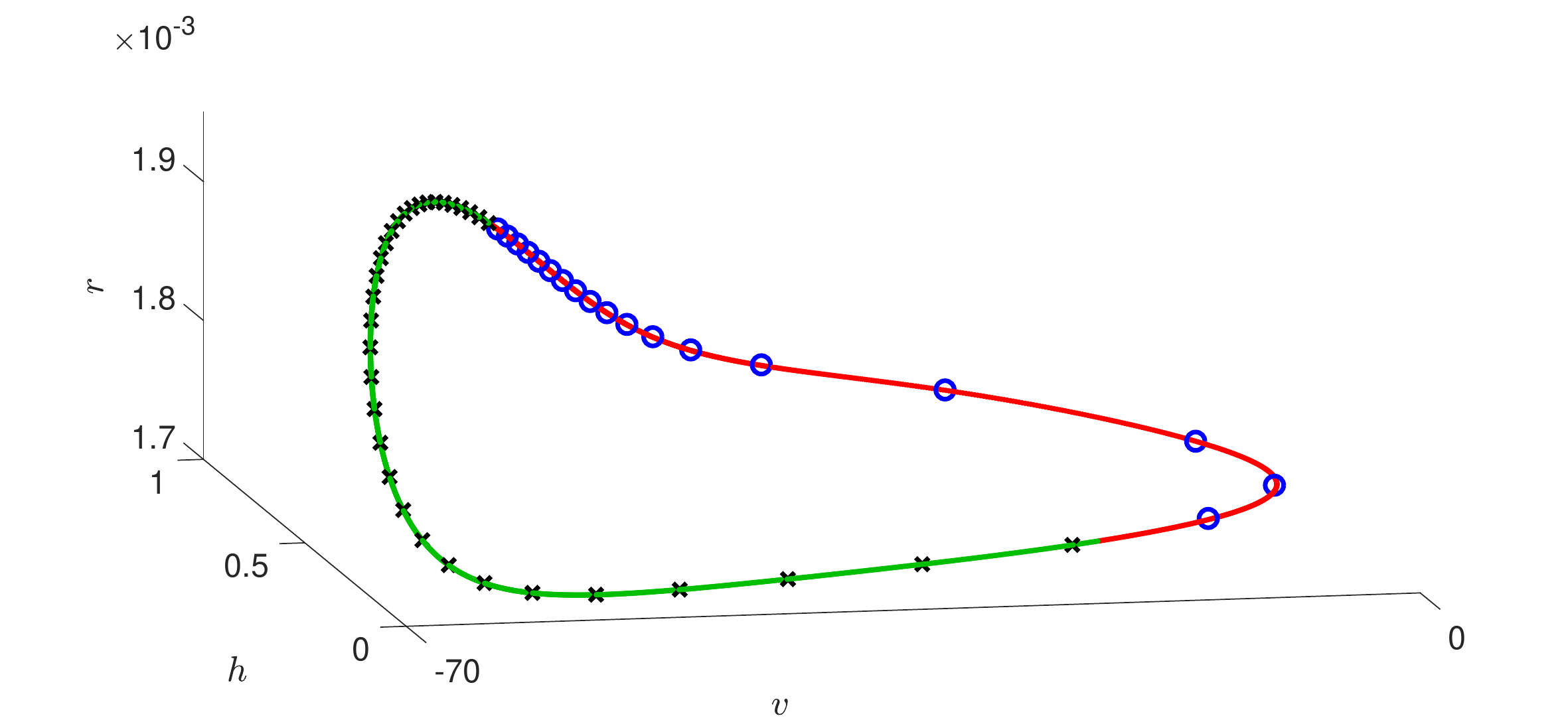}
\end{center}
\caption{Thalamic synchronized population oscillation: The closed curve is the synchronized oscillation. Small blue circles (resp., black $\times$'s) represents elements of the set $\mathcal{X}$ where the control policy given by elements of the set $\mathcal{U}$ is $-u_1$ (resp., $u_1$). The oscillation is plotted in red (resp., green ) where $\mathcal{Z}'(\theta)$ is negative (resp., positive).}
\label{th_pop}
\end{figure}

\begin{figure}[!t]
\begin{center}
\includegraphics[width=\textwidth]{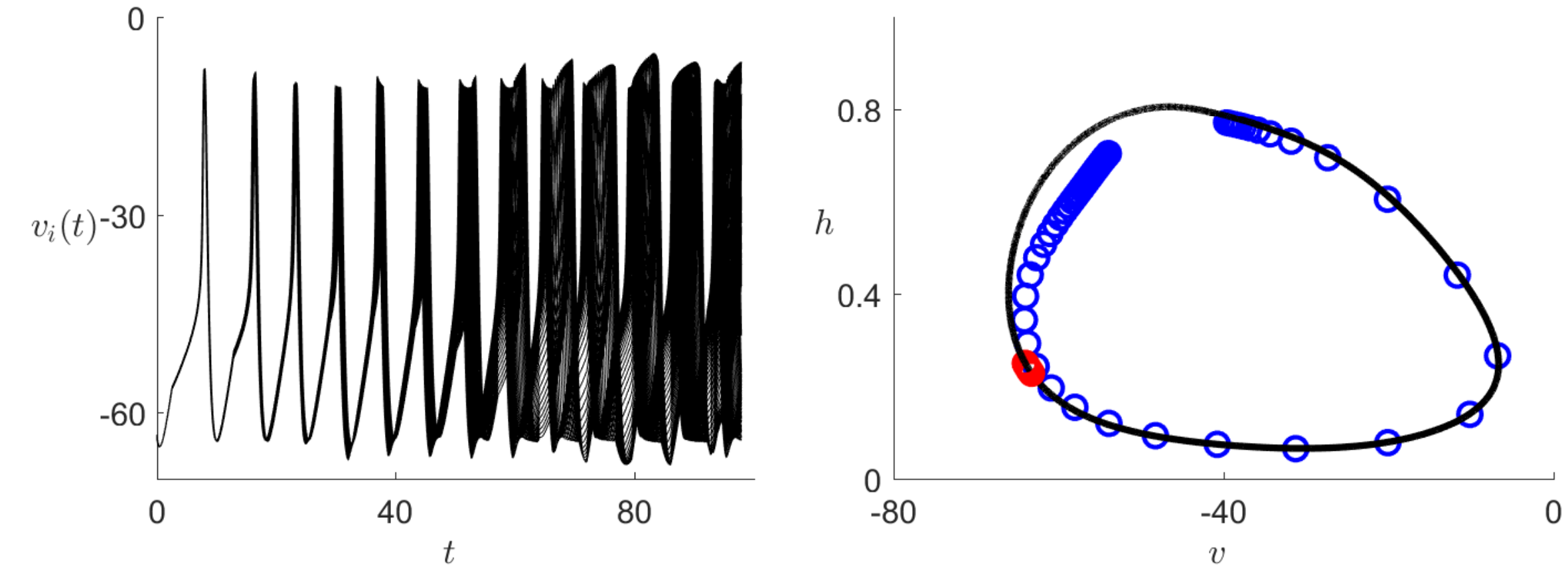}
\end{center}
\caption{Desynchronization of thalamic neuron population: Left panel plots the state $v_i$ for $i=1, \cdots, 51$ neurons as a function of time. Right panel plots the initially synchronized (resp. final desynchronized) neurons as small red (resp., blue) circles.}
\label{pop_result}
\end{figure}

The generated control policy shown in Figure \ref{th_pop} gives a positive control input in the bottom left region of oscillation and a negative control in the top right region of the oscillation.
The same figure also plots a model based control policy discussed below. Figure \ref{pop_result} shows the results of desynchronization of a thalamic neuron population by our learning algorithm. As shown in both the left and right panels of the figure, the control policy from our learning algorithm is able to desynchronize an initially synchronized thalamic neuron population in about 90~ms while keeping the oscillators close to the initially synchronized oscillation. It may seem that the population is largely synchronized from the right panel of Figure \ref{pop_result} but that is not the case. Since the oscillators spend most of their time near the top of the limit cycle, one naturally observes more of them near the top of the limit cycle even though they are evenly spread out in time (and not space). This becomes clear from the left panel of Figure \ref{pop_result}.

\subsection{Model based validation of control policy}
Here we analyze why the policy predicted by our learning algorithm works. Consider an oscillator population comprised of just 2 oscillators whose dynamics evolve according to phase reduction as 
\begin{eqnarray}
\dot \theta_1= \omega + \mathcal{Z}(\theta_1)u(t),\\
\dot \theta_2 = \omega + \mathcal{Z}(\theta_2)u(t).
\end{eqnarray}
The dynamics of their phase difference $\phi = \theta_1 - \theta_2$ can be written as (cf, \cite{dan14})
\begin{equation}
\dot \phi = \mathcal{Z}'(\theta) u(t) \phi +{\cal{O}}(\phi ^3),\label{phase_diff}
\end{equation}
where $\theta=0.5(\theta_1+\theta_2)$ is the mean of the two oscillators' phases, and $\mathcal{Z}'(\theta)$ is the derivative of the phase response curve with respect to $\theta$. If the oscillators are synchronized then their phase difference $\phi \approx 0$, thus higher order term in equation (\ref{phase_diff}) can be ignored and the equation can be rewritten as
\begin{equation}
\dot \phi =\mathcal{Z}'(\theta) u(t) \phi .
\end{equation}
To desynchronize these two synchronized oscillators the coefficient of $\phi$ in the above equation \label{diff_phase} should be positive. This can be achieved if $u(t)$ is of same sign as ${Z}'(\theta)$. This is exactly what our policy predicts, as is shown in Figure \ref{th_pop}. The policy predicts the control to be positive in the region of oscillation where $\mathcal{Z}'(\theta)$ is positive, and it predicts the control to be negative in the region of oscillation where $\mathcal{Z}'(\theta)$ is negative, thus explaining why our algorithm is able to desynchronize the oscillator population.

\section{Stabilizing an Unstable Fixed Point}\label{unstable}

In this section we apply our second learning algorithm to stabilize an unstable fixed point of an underactuated dynamical system. This control objective is one of the oldest studied control theory problems that is employed in several fields including robotics, electrochemical systems, and treatment of cardiac arrhythmias \cite{under_mec3,electrochem,Monga2019}. To demonstrate this, we consider the Lorenz system \cite{lorenz} given as:
\begin{eqnarray}
\dot x&=&\sigma(y-x) +\pi\left(\mbf{x}(t)\right) ,\label{loren1}\\
\dot y &=& rx-y-xz,\\
\dot z &=& xy -bz.
\end{eqnarray}
In the absence of control input with parameters $\sigma=10, \ b=8/3, \ r=1.5$, the system is bistable with $\mbf{x}^{s1}=(-1.15,-1.15,0.5)$, $\mbf{x}^{u}=(0,0,0)$, and $\mbf{x}^{s2}=(1.15,1.15,0.5)$, all shown in Figure \ref{lorenz2}. 

\subsection{Learning Algorithm}
The control objective is for a trajectory to converge to the unstable fixed point $\mbf{x}^{u}$ starting anywhere in the state space. We choose our reward function to be the negative of the Euclidean distance between current state and the desired unstable fixed point
\begin{equation}
\mathcal{R}(\mbf{x}(t)) = -||\mbf{x}(t)-\mbf{x}^{u}||.\label{reward_unstable}
\end{equation}
Thus the control will make the trajectory converge towards the desired fixed point while increasing the reward to 0. To converge to $\mbf{x}^{u}$ starting anywhere in the state space, we use our learning algorithm to generate a control policy. We take $u_1=-u_2=5$, and sample $N=1000$ points for generating the sets $\mathcal{X}, \ \mathcal{U}$. We choose $\Delta t = 0.001$ and take the binary classifier parameter $\tau=5$. 

\subsection{Results}
\begin{figure}[!t]
\begin{center}
\includegraphics[width=\textwidth]{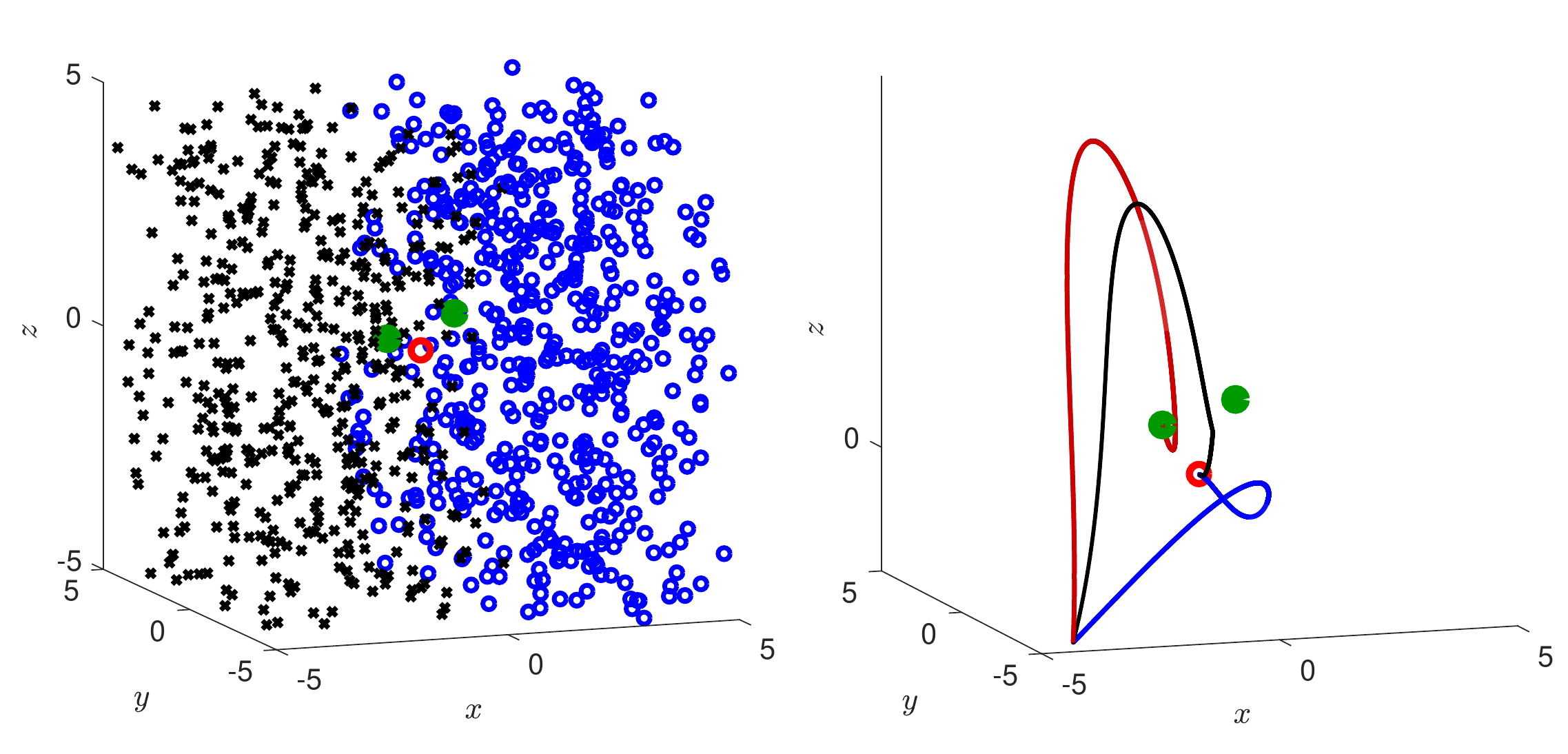}
\end{center}
\caption{Lorenz system: Solid green (resp., open red) circles represent $\mbf{x}^{s1}$, $\mbf{x}^{s2}$, (resp., $\mbf{x}^{u}$). In the left panel, open blue circles (resp., black $\times$'s) represent elements of the set $\mathcal{X}$ where the control policy given by elements of the set $\mathcal{U}$ is -5 (resp., 5). In the right panel, the uncontrolled trajectory plotted in red converges to $\mbf{x}^{s1}$, and the supervised learning (resp., Lyapanov) based control trajectory plotted in black (resp., blue) converges to $\mbf{x}^{u}$.}
\label{lorenz2}
\end{figure}

The generated control policy is shown in the left panel of Figure \ref{lorenz2}. 
Blue open circles (resp., black $\times$'s) represent elements of the set $\mathcal{X}$ where the control policy given by elements of the set $\mathcal{U}$ is -5 (resp., 5). A controlled trajectory using our learning algorithm and an uncontrolled trajectory starting from same $\mbf{x}(0)$ are shown in the right panel of Figure \ref{lorenz2}. The learning algorithm converges the trajectory to $\mbf{x}^{u}$, whereas the uncontrolled trajectory converges to $\mbf{x}^{s1}$. In doing so, the learning based control consumes 150 units of control energy ($\int_0^6 \pi\left(\mbf{x}(t)\right)^2 dt$). We investigate the robustness of our learning algorithm by testing it on 1000 randomly generated initial conditions, and in all the 1000 cases, the learning based control algorithm is able to converge the trajectories within a ball of radius 0.09 units centered at $\mbf{x}^{u}$, achieving 100\% effectiveness. 

\subsection{Comparison with Lyapanov based control}
To demonstrate energy efficiency of our learning algorithm, we compare it with Lyapunov-based control to stabilize $\mbf{x}^{u}$. Consider the following positive definite Lyapunov function
\begin{equation}
V(t) = \frac{1}{2}x(t)^2 + \frac{1}{2}y(t)^2 + \frac{1}{2}z(t)^2.
\end{equation}
Its time derivative is given as 
\begin{equation}
\dot V(t) = -2\sigma x(t)^2 -2y(t)^2 -2bz(t)^2 +2x(t)\left(u(t)+(\sigma + r)y(t)\right),
\end{equation}
where $u(t)$ takes the place of $\pi\left(\mbf{x}(t)\right)$ in equation (\ref{loren1}). Then by taking $u(t) = -(\sigma + r)y(t)$, one gets a negative definite time derivative of the Lyapunov function. Thus by the Lyapunov theorem, this control asymptotically stabilizes the unstable fixed point $\mbf{x}^{u}$. The control trajectory based on this control is plotted in blue in the right panel of Figure \ref{lorenz2}. As seen in the figure, the Lyapunov-based control is able to converge the trajectory towards the unstable fixed point as well. But in doing so, it consumes 1176.8 units of control energy ($\int_0^6 u(t)^2 dt$), which is almost 8 times the energy consumed by our learning based control. This is partly because our learning based control uses the inherent system dynamics to control the trajectory, as the controlled trajectory seems to stay close to the uncontrolled trajectory. In contrast, the Lyapanov based control drives the trajectory far away before it converges to $\mbf{x}^{u}$, thus it ends up consuming much more energy.

\section{Robustness to Noise}\label{noise}
We have demonstrated the effectiveness of our algorithms in several scenarios in which the algorithms were based on data generated from a deterministic dynamical model. However, real data measured from an experimental setup will be noisy. In order for our algorithms to work in an experimental setup it is imperative to investigate their performance when the data is corrupted with noise. We do that by considering the Duffing oscillator in the bistable parameter regime from Section \ref{duff_no_noise}. The control objective is still for the trajectory to converge to $\mbf{x}^{s2}$ starting anywhere in the state space.

\subsection{Learning Algorithm}
\begin{figure}[!t]
\begin{center}
\includegraphics[width=0.6\textwidth]{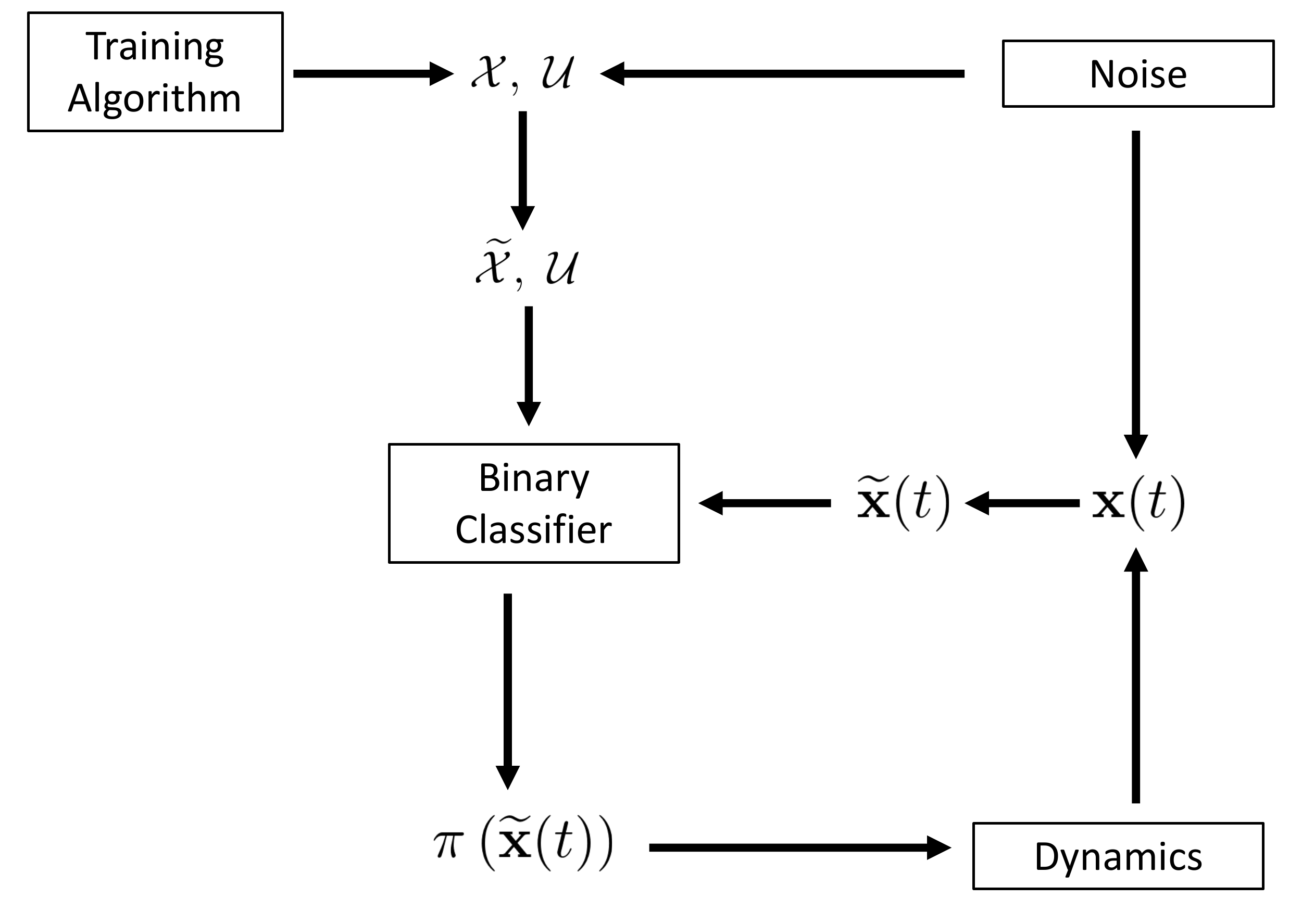}
\end{center}
\caption{Flowchart of the Supervised Learning Algorithm with added noise.}
\label{flowchart1}
\end{figure}

To replicate the effect of noise in an experimental setup, we use exactly the same parameters as before to generate the sets $\mathcal{X}, \ \mathcal{U}$ and corrupt the set $\mathcal{X}$ by adding Gaussian white noise with mean 0 and standard deviation $\sigma$, resulting in the set $\widetilde{\mathcal{X}}$. Thus each element in the dataset will be offset from its true value. We also add Gaussian white noise with the same mean and standard deviation to the state $\mbf{x}(t)$ resulting in $\widetilde{\mbf{x}}(t)$, which the binary classifier takes as input at every time step. This accounts for the noise in estimation of the state by the classifier in a real system. The flowchart from Figure \ref{flowchart} with added noise is modified and shown in Figure \ref{flowchart1}.

\subsection{Results}
\begin{figure}[!t]
\begin{center}
\includegraphics[width=\textwidth]{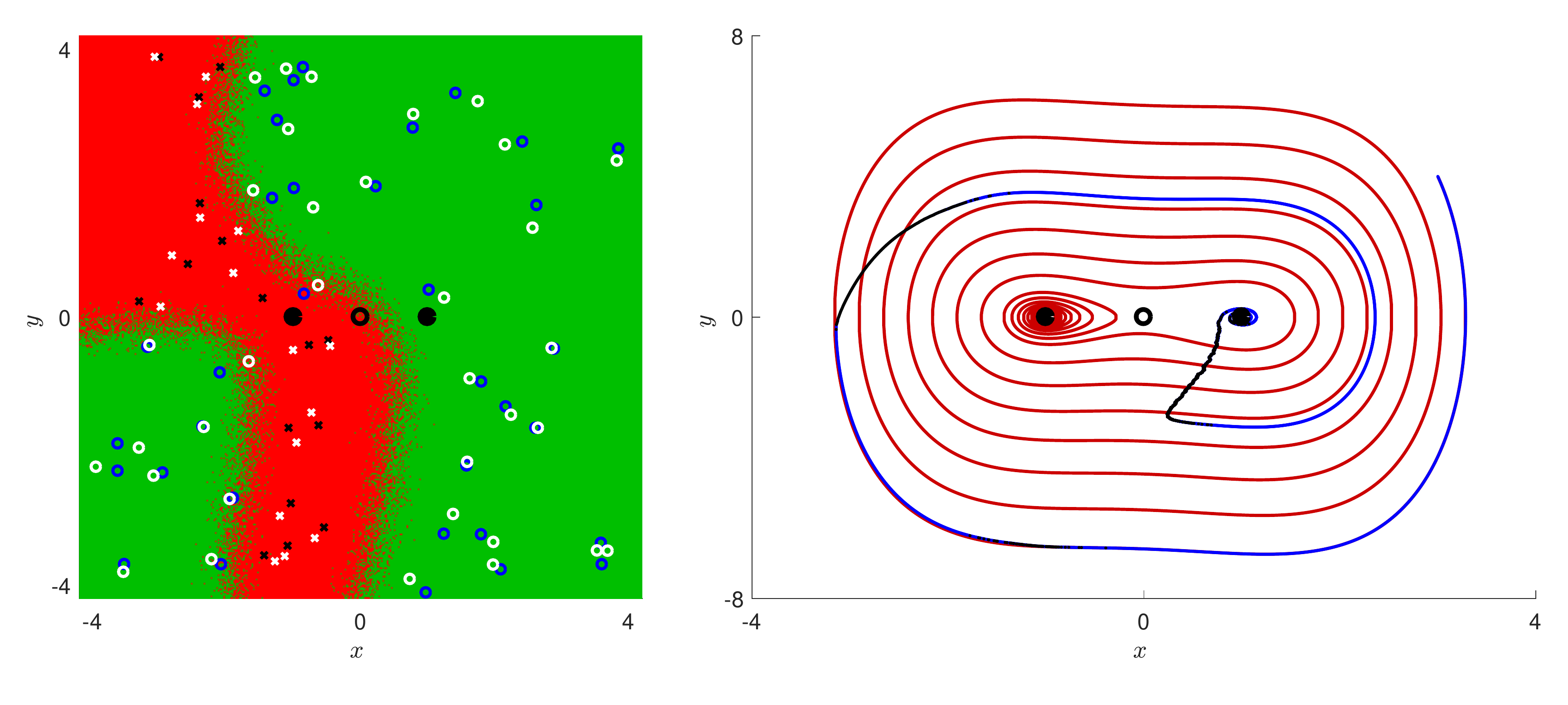}
\end{center}
\caption{Duffing Oscillator with noise ($\delta=0.1$): Solid (resp., open) black circles represent $\mbf{x}^{s1}$, $\mbf{x}^{s2}$, (resp., $\mbf{x}^{u}$). In the left panel, open blue circles (resp., black $\times$'s) represent elements of the set $\widetilde{\mathcal{X}}$ where the control policy given by elements of the set $\mathcal{U}$ is OFF (resp., ON). The green (resp., red) region is where the output $\pi\left(\widetilde{\mbf{x}}(t)\right)$ of the binary classifier is OFF (resp., ON). The elements of the original set $\mathcal{X}$ are plotted in white. In the right panel, the trajectory starts in the region of attraction of $\mbf{x}^{s1}$, and converges to $\mbf{x}^{s2}$ (resp., $\mbf{x}^{s1}$) with (resp., without) control. When the control policy is ON (resp., OFF), the trajectory is plotted in black (resp., blue). The uncontrolled trajectory is plotted in red.}
\label{duff1}
\end{figure}

The generated control policy corrupted with noise of standard deviation $\sigma=0.2$  is shown in the left panel of Figure \ref{duff1}. 
Blue open circles (resp., black $\times$'s) represent elements of the noise corrupted set $\widetilde{\mathcal{X}}$ where the control policy given by elements of the set $\mathcal{U}$ is OFF (resp., ON). The green (resp., red) region is where output $\pi\left(\widetilde{\mbf{x}}(t)\right)$ of the binary classifier is OFF (resp., ON). The elements of the original set $\mathcal{X}$ are plotted in white to show the shifting of the elements due to the noise. Besides shifting the elements, the addition of white noise blurs the decision boundary between the ON and OFF policy region . A controlled and an uncontrolled trajectory starting from the same $\mbf{x}(0)$ are shown in the right panel of Figure \ref{duff1}. As can be seen in this figure, the control algorithm converges the trajectory to $\mbf{x}^{s2}$ even though it has been corrupted by adding noise both to the training dataset and to the input of of the binary classifier. On the other hand, the uncontrolled trajectory converges to $\mbf{x}^{s1}$. Note that the controlled trajectory follows a different route when compared to the noiseless case from Figure \ref{duff}. This is because the noise distorts the decision boundary between ON and OFF states of the policy, resulting in a slightly different path.

\subsection{Robustness and Noise Intensity}
We first investigate the robustness of our learning algorithm in the presence of noise by considering noise with intensity $\sigma = 0.2$ and testing it on 1000 randomly generated initial conditions; in all 1000 cases the control algorithm is able to converge the trajectories to $\mbf{x}^{s2}$, achieving 100\% effectiveness. It is natural to ask  if effectiveness decreases down from 100\% as one increases the noise intensity beyond 0.2. To answer this question, we tested robustness of our algorithm by simulating 1000 randomly generated initial conditions and evolving them under our control policy for various levels of noise intensity $\sigma$ and bandwidth parameter $\tau$. We found that our algorithm can still achieve very high effectiveness for higher $\sigma$, if we adjust the bandwidth parameter $\tau$.
    
\begin{figure}[!t]
\begin{center}
\includegraphics[width=\textwidth]{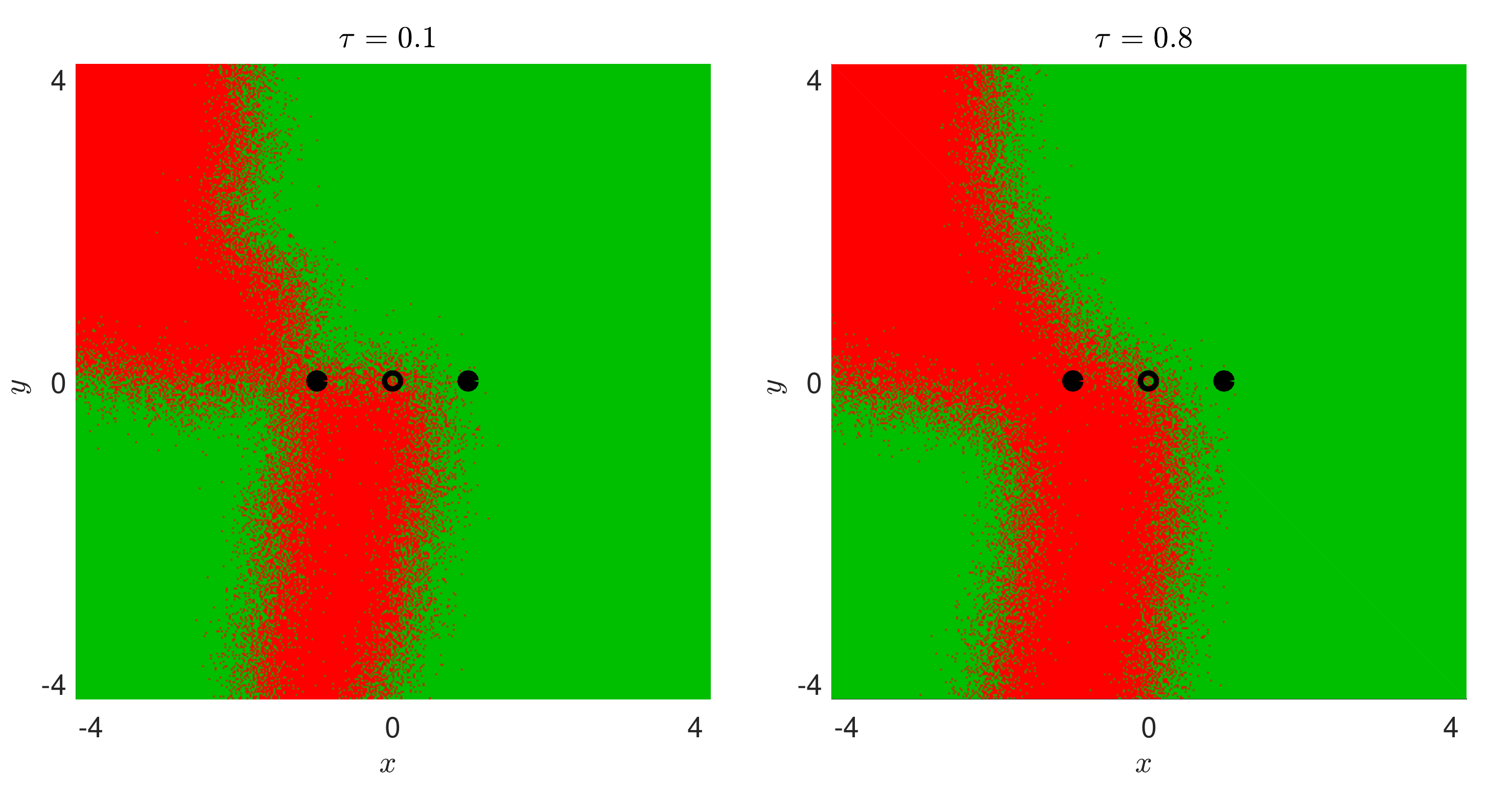}
\end{center}
\caption{Duffing Oscillator policy with noise intensity $\sigma=0.3$.}
\label{policy_tau}
\end{figure}

To understand why this is the case, see Figure \ref{policy_tau},
where we plot two policies for noise intensity $\sigma=0.3$ with bandwidth parameter $\tau = 0.1$ and $0.8$. Since noise randomly shifts elements of our data set $\mathcal{X}$, it distorts the decision boundary when compared to the decision boundary boundary in the noiseless case (see the left panel of Figure \ref{duff} for a comparison). For small $\tau$ values, this leads to significant segmentation of the policy  as the output of the binary classifier can vary rapidly with a small change in the input state. This is evident from the left panel of Figure \ref{policy_tau}. On the other hand, if we increase the bandwidth parameter, segmentation is significantly reduced, as is evident from the right panel of the same figure. A higher bandwidth parameter leads to less segmentation, because the binary classifier effectively takes into account more neighboring data points in making its decision (equation (\ref{bc1})). In other words, when these points are randomly shifted by the addition of Gaussian white noise, the weighted average from a larger number of data points produces smoother transition between regions, and thus leads to reduced segmentation. On the other hand, if the bandwidth parameter was too small, the binary classifier's decision would be based on only a very few neighboring data points, which, when randomly shifted by Gaussian white noise, can produce a highly segmented and non-smooth policy.

But robustness does not keep increasing with increasing $\tau$. A $\tau$ value which is too high can actually reduce effectiveness, as depicted in Table \ref{robust_tau}, because, as previously mentioned, high $\tau$ values can wrongly extend the decision boundary in favor of a cluster of data points which is too large. Thus, a carefully chosen moderate value of $\tau$ is most effective in assuring robustness in the presence of noise. Another observation from our robustness studies for various levels of noise intensity $\sigma$ and $\tau$ depicted in Table \ref{robust_tau} is that we need an increasing value of $\tau$ to achieve 100\% effectiveness as we increase the noise intensity $\sigma$. This is also intuitive to understand, as for a given value of $\tau$, one expects more segmentation in the control policy with increasing $\sigma$. Thus, the locally weighted nature of our binary classifier, together with the flexibility to adjust its bandwidth allows our algorithm to be effective in the presence of noise.

\begin{table}[!t]
 \centering 
\caption{Effectiveness of our algorithm as a function of noise intensity $\sigma$ and bandwidth parameter $\tau$ for the Duffing Oscillator.}
\label{robust_tau}
\renewcommand{\arraystretch}{1.5}

\begin{tabular}[c]{c | ccccc}
\tikz[diag text/.style={inner sep=2pt, font=\large},
      shorten/.style={shorten <=#1,shorten >=#1}]{%
        \node[below left, diag text] (sigma) {$\sigma$};
        \node[above right=2pt and -2pt, diag text] (tau) {$\tau$};
        \draw[shorten=4pt, thin] (sigma.north west|-tau.north west) -- (sigma.south east-|tau.south east);}
 & 0.1 & 0.4 & 0.8 & 1.2 & 1.6\\
\hline
0.2 & 100\% & 100\% & 100\%  & 100\%  & 100\% \\
0.3 & 100\% & 100\% & 100\% & 100\% & 83\% \\ 
0.4 & 0\% & 89.9\%  & 100\% & 93.2\% & 64.1\%\\ 
0.5 & 0\% & 78.6\% & 90.8\% &100\% & 68.7\% \\
0.6 & 0\%  & 0\%  & 0\%  & 9\% & 95.6\% \\ 
\end{tabular}
\end{table}

While it may not be possible to achieve 100\% effectiveness in all control applications, we believe that this example illustrates key considerations in characterizing the robustness of our algorithms to noise. Another possible consideration for achieving high effectiveness could be choosing a higher sample size $N$, however this would increase the computational cost of our algorithm. Since we are able to achieve very high effectiveness by adjusting the bandwidth parameter while keeping the computational cost the same, we favor this over adjusting $N$ in characterizing the robustness of our algorithms to noise.

\section{Conclusion}\label{conclusion}
In this article we have developed two novel supervised learning algorithms to control a range of underactuated dynamical systems. The algorithms output a bang bang (binary) control input to achieve the desired control objectives which maximize a reward function. A simple yet intelligent structure allows the algorithms to be energy efficient as they learn to take advantage of the inherent dynamics. We demonstrated the versatility of our algorithms by applying them to a diverse range of applications including: switching between bistable states, changing the phase of an oscillator, desynchronizing a population of synchronized coupled oscillators, and stabilizing an unstable fixed point. For most of these applications we were able to reason why our algorithms work by using traditional dynamical systems and control theory. We also compared our algorithms to some traditional nonlinear model control algorithms and showed that our algorithms work better. We also carried out a robustness study to demonstrate the effectiveness of our algorithms even with noisy data.

We simulated various dynamical models to generate data for training our supervised learning algorithms. In an experimental setting such an algorithm can be similarly implemented by stimulating the system with binary control inputs at different states of the system and determining which control input works best for the different sampled states, and ultimately using that information in constructing a binary classifier. The data generated from an experimental setting might be corrupted with noise, thus to demonstrate the potential of our algorithm in a real setting we showed that our algorithm works even in the presence of noise. The unique structure of our binary classifier comes to the rescue by negating the adverse effects of noise. Note that having an additive control input doesn't restrict our algorithms. Since the algorithms work by maximizing the reward function, the structure of the control input coming into the dynamics does not matter. In the future, we plan to explore how to modify our algorithm if some of the states are not observable, and how to adapt our algorithm for very high dimensional dynamical systems.

\appendix
\section{Models}\label{model_param}
In this appendix, we give details of the mathematical models used in this article.

\subsection{Reduced Hodgkin-Huxley model}\label{hh_param}
Here we give the reduced Hodgkin-Huxley model \cite{canard,Keener2009,huxley} used in Section \ref{hh}:
\begin{eqnarray*}
\dot v&=& \left(I - g_{Na}(m_\infty(v))^3(0.8-n)(v-v_{Na}) - g_Kn^4(v-v_K)-g_L(v-v_L)\right)/c+\pi\left(\mbf{x}(t)\right),\\ 
\dot n &=& a_n(v)(1-n)-b_n(v) n,
\end{eqnarray*}
where $v$ is the trans-membrane voltage, and $n$ is the gating variable. $I$ is the baseline current, which we take as 6.69 $\mu A/cm^2$, and $\pi\left(\mbf{x}(t)\right)$ represents the applied control current.
\begin{eqnarray*}
a_n(v) &=&0.01(v+55)/(1-\exp(-(v+55)/10)),\\ 
b_n(v) &=& 0.125\exp(-(v+65)/80),\\
a_m(v) &=& 0.1(v+40)/(1-\exp(-(v+40)/10)),\\
b_m(v) &=& 4\exp(-(v+65)/18),\\
m_\infty(v) &=& a_m(v) / (a_m(v) + b_m(v)),\\
c &=& 1,\   g_L = 0.3,\      g_{Na} = 120,\    v_{Na} = 50\\g_K &=36& ,\    v_K = -77,\     v_L = -54.4\ I =20 .
\end{eqnarray*}

\subsection{Thalamic neuron model}\label{thalam_param}
The thalamic neuron model is given as
\begin{eqnarray*}
\dot v&=&\frac{-I_L(v)-I_{Na}(v,h)-I_K(v,h)-I_T(v,r)+I_b }{C_m}+\pi\left(\mbf{x}(t)\right),\\
\dot h&=&\frac{h_{\infty}(v)-h}{\tau_h(v)},\\
\dot r&=&\frac{r_{\infty}(v)-r}{\tau_r (v)}.
\end{eqnarray*}
where
\begin{eqnarray*}
    h_\infty (v)&=& 1/(1+\exp((v+41)/4)),\\
    r_\infty (v)&=& 1/(1+\exp((v+84)/4)),\\
    \alpha_h(v) &=& 0.128\exp(-(v+46)/18),\\
    \beta_h (v)&=& 4/(1+\exp(-(v+23)/5)),\\
    \tau_h(v)&=& 1/(\alpha_h+\beta_h),\\
    \tau_r (v)&=& (28+\exp(-(v+25)/10.5)),\\
    m_\infty (v)&=& 1/(1+\exp(-(v+37)/7)),\\
    p_\infty(v) &=& 1/(1+\exp(-(v+60)/6.2)),
    \end{eqnarray*}
    \begin{eqnarray*}
    I_L(v)&=&g_L(v-e_L),\\
    I_{Na}(v,h)&=&g_{Na}({m_\infty}^3)h(v-e_{Na}),\\
    I_K(v,h)&=&g_K((0.75(1-h))^4)(v-e_K),\\
    I_T(v,r)&=&g_T(p_\infty^2)r(v-e_T),
\end{eqnarray*}
\[
C_m = 1, \qquad   g_L = 0.05, \qquad   e_L = -70, \qquad    g_{Na} = 3, \qquad   e_{Na} = 50,
\]    
\[
g_K = 5,\qquad    e_K = -90,\qquad    g_T = 5,\qquad    e_T = 0, \qquad I_b = 5.
\]

\section{Phase Reduction}\label{phase_red}
Phase reduction is a classical technique to describe dynamics near a limit cycle. It works by reducing the dimensionality of a dynamical system to a single phase variable $\theta$~\cite{winfree,kuramoto}. Consider a general $n$-dimensional dynamical system given by 
\begin{equation*}
\frac{d \mbf{x}(t)}{dt} = F(\mbf{x}(t)), \qquad \mbf{x}(t) \in \mathbb{R}^n, \qquad (n \ge 2) .
\label{dxdt}
\end{equation*}
Suppose this system has a stable periodic orbit $\mbf{x}^{s}(t)$ with period $T$. For each point $\mbf{x}^*$ in the basin of attraction of the periodic orbit, there exists a corresponding phase $\theta(\mbf{x}^*)$ such that
\begin{equation*}
\lim_{t \rightarrow \infty} \left| \mbf{x}(t) - \mbf{x}^{s} \left( t+\frac{T}{2 \pi} \; \theta(\mbf{x}^*) \right) \right| = 0,
\end{equation*}
where $\mbf{x}(t)$ is the flow of the initial point $\mbf{x}^*$ under the given vector field. The function $\theta(\mbf{x})$ is called the {\it asymptotic phase} of $\mbf{x}$, and takes values in 
$[0, 2 \pi)$.  For neuroscience applications, we typically take $\theta=0$ to correspond to the neuron firing an action potential.  {\it Isochrons} are level sets of this phase function, and it is typical to define isochrons so that the phase of a trajectory advances linearly in time both on and off the limit cycle, which implies that
\begin{equation*}
\frac{d\theta}{dt} = \frac{2 \pi}{T} \equiv \omega
\label{theta_flow}
\end{equation*}
in the entire basin of attraction of the limit cycle. Now consider our underactuated system given by equation (\ref{udxdt}). Phase reduction can be used to reduce this system to a one-dimensional system given by \cite{winf01,kura84,brow04,tutorial}:
\begin{equation*}
\dot \theta = \omega + \mathcal{Z}(\theta)\pi\left(\mbf{x}(t)\right).
\end{equation*}
Here $\mathcal{Z}(\theta)$ is the first component of the gradient of phase variable $\theta$ evaluated on the periodic orbit and is referred to as the {\it (infinitesimal) phase response curve (PRC)}. It quantifies the effect of an infinitesimal control input on the phase of a limit cycle.

\section*{Acknowledgment}
This work was supported by National Science Foundation Grant No. NSF-1635542.

\bibliographystyle{siamplain}
\bibliography{ms}

\end{document}